\documentclass[12pt]{article}
\usepackage{amssymb}
\newcommand{\sect}[1]{\setcounter{equation}{0}\section{#1}}

\newcommand{\be}{\begin{equation}}
\newcommand{\ee}{\end{equation}}
\newcommand{\bea}{\begin{eqnarray}}
\newcommand{\eea}{\end{eqnarray}}
\newcommand{\beano}{\begin{eqnarray*}}
\newcommand{\eeano}{\end{eqnarray*}}
\newcommand{\nonu}{\nonumber \\}

\newcommand{\ca}{\mbox{$\cal{A}$}}
\newcommand{\cb}{\mbox{$\cal{B}$}}

\newcommand{\cf}{\mbox{${\cal F}$}}
\newcommand{{\cg}}{\mbox{$\cal{G}$}}

\newcommand{\cR}{\mbox{$\cal{R}$}}
\newcommand{\cs}{\mbox{$\cal{S}$}}

\newcommand{\cu}{\mbox{${\cal U}$}}
\newcommand{\cv}{\mbox{${\cal V}$}}

\newcommand{\wh}[1]{\widehat{#1}}
\newcommand{\wt}[1]{\widetilde{#1}}
\newcommand{\mb}[1]{\ \mbox{ #1 }\ }
\newcommand{\half}{\frac{1}{2}}

\newtheorem{prop}{Property}[section]
\newtheorem{coro}[prop]{Corollary}
\newtheorem{defi}[prop]{Definition}
\newtheorem{theo}[prop]{Theorem}

\newtheorem{lem}[prop]{Lemma}
\newcommand{\prf}{\underline{Proof:}\ }
\newcommand{\finprf}{\null \hfill {\rule{5pt}{5pt}}\\[2.1ex]\indent}
\newcommand{\ie}{{\it i.e.}\ }
\newcommand{\CC}{\mbox{${\mathbb C}$}}

\newcommand{\II}{\mbox{${\mathbb I}$}}

\newcommand{\CMP}[1]{Commun.\ Math.\ Phys.\ {\bf #1}}

\begin{document}
\renewcommand{\thefootnote}{\fnsymbol{footnote}}
\newpage
\pagestyle{empty}
\setcounter{page}{0}

\newcommand{\LAP}{LAPTH}
\def\logo{{\bf {\huge LAPTH}}}

\centerline{\logo}

\vspace {.3cm}

\centerline{{\bf{\it\Large 
Laboratoire d'Annecy-le-Vieux de Physique Th\'eorique}}}

\centerline{\rule{12cm}{.42mm}}

\vspace{20mm}

\begin{center}
    
{\LARGE  {\sffamily Quantum group symmetry of integrable 
 systems\\[1.2ex] 
 with or without boundary }}\\[1cm]

{\large E. Ragoucy\footnote{ragoucy@lapp.in2p3.fr}}\\[.21cm] 
 Laboratoire de Physique Th{\'e}orique \LAP\footnote{UMR 5108 
 du CNRS, associ{\'e}e {\`a} l'Universit{\'e} de Savoie.}\\[.242cm]
 LAPP, BP 110, F-74941  Annecy-le-Vieux Cedex, France. 
\end{center}
\vfill\vfill

\begin{abstract}
We present a construction of integrable hierarchies without or with 
boundary, starting from a single $R$-matrix, or equivalently from a 
ZF algebra. We give explicit 
expressions for the Hamiltonians and the integrals of motion of the 
hierarchy in term of the ZF algebra. In the case without boundary, 
the integrals of motion form a quantum group, while in the case with 
boundary they form a Hopf coideal subalgebra of the quantum group.
\end{abstract}

\vfill
\centerline{MSC-classification: 81R10, 81R12}
\vfill
\begin{center}
This paper is based on talks given at

-- 5$^{th}$ Bologna Workshop on {\it Conformal Field Theories 
and Integrable Models}, 
Bologna (Italy), Sept. 2001

-- LAQ'2001, {\it Australian Research Symposium on Groups and their 
Representations}, Auckland (New Zealand), Dec. 2001
\end{center}
\vfill
\rightline{\tt mathQA/0202095}
\rightline{\LAP-900/02}
\rightline{February 02}

\newpage
\pagestyle{plain}

\sect{Introduction}
The aim of the present article is to show that, in the study of 
integrable systems (without or with boundary), all the relevant 
informations are contained in (and can be reconstructed from) a 
unique algebra: the ZF algebra. This fact is known in the case of
the Non Linear Schr\"odinger equation in 2 
dimensions, and we will show that this is true in general.

\null

In the study of the Non Linear Schr\"odinger (NLS) equation in 2 
dimensions, it has 
been known for a long time that a central role is played by a 
deformed oscillators algebra, the ZF algebra \cite{ZF}. These deformed
oscillators can be seen as the asymptotic states of the physical 
system. Thanks to this algebra, 
one can reconstruct the quantum Hamiltonians of the hierarchy 
associated to NLS, as well as its quantum canonical fields \cite{qNLS}. Quite 
recently \cite{MRSZ}, it has also been shown that 
the integrals of motion of this 
hierarchy, which form a Yangian algebra $Y(N)$, are a part of the ZF 
algebra, and an explicit construction in term of basic ZF generators 
was given. 

It appears that the same technics can be applied in the case of the 
Non Linear Schr\"odinger equation with a boundary (BNLS). Following the
work of Cherednik \cite{Cher} and Sklyanin \cite{Skly}, who studied
the problem of boundaries in integrable systems in the QISM framework, 
Mintchev et al \cite{Mint,bound} have shown that the algebraic approach
used for NLS can also be applied to  BNLS. In that case, the ZF algebra
is replaced by a so-called boundary algebra, which contains both the
asymptotic states of the system and the effect of the boundary. 
Indeed, the boundary algebra contains as a subalgebra the reflection 
algebra, which is known to be of fundamental importance in 
integrable  systems with boundary. The Hamiltonians  and the canonical 
fields of the hierarchy
can be also reconstructed from the boundary algebra, and recently 
\cite{BNLS}  it 
has been shown that the reflection algebra indeed correspond to the 
integrals of motion of the hierarchy. Thanks to this subalgebra, a 
classification of the boundary conditions was also given \cite{BNLS}.

\null

We will show that most of the above features, explicited for NLS, are 
valid in full generality, provided one has at its disposal 
a $R$-matrix which obeys a 
unitary condition. In particular, the Hamiltonians and the 
integrals of motion of an integrable hierarchy will be exhibited, both 
in the case without and with boundary. 
We will show that the integrals 
of motion form a quantum group when there is no boundary, and a Hopf 
coideal subalgebra of this quantum group when there is a boundary. 
As for NLS, the central role of the construction  will be a ZF 
algebra, and a class of operators (contained in the ZF algebra) 
called well-bred operators.

The paper is organized as follows: in the first section, we treat the 
case without boundary; then, the next section deals with the case of 
boundaries; finally, we conclude in section \ref{conc}. Most of the 
results presented here can be found in \cite{ZFqgp,ZFbound}, 
and we refer to these 
papers for the original proofs.

\sect{Case without boundary}
The starting point is an evaluated $R$-matrix, of size $N^2\times 
N^2$, with spectral parameter, and which obeys the Yang-Baxter equation 
and the unitarity condition:
\bea
&& R_{12}(k_1,k_2)R_{13}(k_1,k_3)R_{23}(k_2,k_3)=
R_{23}(k_2,k_3)R_{13}(k_1,k_3)R_{12}(k_1,k_2)\label{YBE}
\\
&& R_{12}(k_1,k_2)R_{21}(k_2,k_1)=\II\otimes \II \label{unitarity}
\eea
Note that no assumption is made on the form of the dependence in the 
spectral parameters $k_{1}$ and $k_{2}$.

To this $R$-matrix one can introduce two types of algebras, that we 
describe below.
\subsection{The Z.F. algebra}
The Zamolodchikov-Faddeev (ZF) algebra is an exchange algebra, that can 
be seen as a deformation of an oscillators algebras.
\begin{defi}[ZF algebra $\ca_{R}$]\hfill\\
To the above $R$-matrix, one can associate a ZF algebra $\ca_{R}$, with generators 
$a_{i}(k)$ and $a^\dag_{i}(k)$ ($i=1,\ldots,N$) and exchange relations:
\bea
a_{1}(k_1)\, a_{2}(k_2) &=& R_{21}(k_{2},k_{1}) \, a_{2}(k_2)\,a_{1}(k_1) 
\label{AN-1}\\
a_{1}^\dag(k_1)\,a_{2}^\dag (k_2) &=& 
a_{2}^\dag (k_2)\,a_{1}^\dag (k_1)\,R_{21}(k_{2},k_{1})  \\
a_{1}(k_1)\,a_{2}^\dag (k_2) &=& a_{2}^\dag (k_2)\,
R_{12}(k_{1},k_{2})\,a_{1}(k_1) 
+\delta(k_{1}-k_{2})\delta_{12}
\label{AN-3}
\eea
\end{defi}
We use the notations on auxiliary spaces
\bea
&&a_{1}(k_1)=\sum_{i=1}^Na_{i}(k_{1})\, e_{i}\otimes\II
\mb{,} a_{2}(k_2)=\sum_{i=1}^Na_{i}(k_{2})\, \II\otimes e_{i}\label{not1}\\
&&a^\dag_{1}(k_{1})=\sum_{i=1}^Na^\dag_{i}(k_{1})\, e^\dag_{i}\otimes\II
\mb{,} a^\dag_{2}(k_{2})=\sum_{i=1}^Na^\dag_{i}(k_{2})\, 
\II\otimes e^\dag_{i}\\
&& \delta_{12}=\sum_{i=1}^N\, e_{i}\otimes e^\dag_{i} 
\mb{,} e^\dag_{i}=(0,\ldots,0,\stackrel{i}{1},0,\ldots,0)
\mb{,} e^\dag_{i}\cdot e_{j}=\delta_{ij}\ \ \label{not3}
\eea
where $\cdot$ stands for the scalar product of vectors.
Note that any central element of $\ca_{R}$ is a constant, and that 
this algebra admits an adjoint operation:
\begin{prop}
Let $\dag$ be defined by
\be
\dag\left\{\begin{array}{lcl}
\ca_{R} & \rightarrow & \ca_{R}\\[.23ex]
{a}(k) & \mapsto & {a}^\dag(k)\\[.23ex]
{a}^\dag(k) & \mapsto & {a}(k)
\end{array}\right.\label{dag}
\ee
together with $(xy)^\dag=y^\dag x^\dag$, $\forall\, 
x,y\in\ca_{R}$ and 
$R_{12}(k_{1},k_{2})^\dag=R_{21}(k_{2},k_{1})$. 

Then $\dag$ is an anti-automorphism of $\ca_{R}$.
\end{prop}
Note that the adjoint operation is defined on $a(k)$ and $a^\dag(k)$, 
\ie on elements of $\ca_{R}[[k]]\otimes \CC^{N}$: that is the reason 
why it has also to be defined on the evaluated $R$-matrix.

\subsection{Quantum group $\cu_{R}$}
Still with the $R$-matrix at our disposable, we can construct a 
quantum group $\cu_{R}$, with Hopf structure, in the usual way. 
\begin{defi}[Quantum group $\cu_{R}$]\hfill\\
To the above $R$-matrix, one can associate a quantum group $\cu_{R}$, 
with generators $T_{ij}^{(n)}$ which  are gathered (using a spectral 
parameter $k$)
in a $N\times N$ matrix 
\be
T(k)=\sum_{i,j=1}^N\sum_{n=0}^\infty k^{-n}
T_{ij}^{(n)}\,E_{ij}
\ee
 which is submitted to the relation
\be
R_{12}(k_{1},k_{2})T_{1}(k_{1})T_{2}(k_{2})=
T_{2}(k_{2})T_{1}(k_{1})R_{12}(k_{1},k_{2})
\ee
with $T_{1}(k)=T(k)\otimes\II$ and $T_{2}(k)=\II\otimes T(k)$.

The coproduct in $\cu_{R}$ is given by
\be
\Delta T(k)=T(k)\stackrel{.}{\otimes} T(k) \label{deltaT}
\ee
where $\stackrel{.}{\otimes}$ is the tensor product in $\cu_{R}$ and 
the matricial product in the auxiliary space.
\end{defi}
More explicitly, the coproduct reads
\be
\Delta T^{ab}_{(n)}=\sum_{p+q=n}\sum_{c=1}^{N}
T^{ac}_{(p)}\otimes T^{cb}_{(q)}
\ee

Depending on the $R$-matrix, one gets for $\cu_{R}$, e.g. the Yangian 
$Y(N)$ based on $gl(N)$ (for 
$R(k_{1},k_{2})=R(k_{1}-k_{2})$ with 
$R(u)=\frac{u^{2}+g^{2}}{u^{2}}(\II-\frac{ig}{u}P_{12})$,
where $P_{12}$ is the permutation of the two auxiliary spaces), 
or the quantum group 
$\cu_{q}(\hat{gl_{N}})$, based on the affine algebra $\hat{gl_{N}}$
(for $R(k_{1},k_{2})=R(k_{1}/k_{2})$ as given in e.g. \cite{glNaff}).

\subsection{Vertex operator construction}
Our aim is to construct $\cu_{R}$ from $\ca_{R}$, \ie $T(k)$ from 
$a(k)$ and $a^\dag(k)$. The basic notion for such a purpose is:
\begin{defi}[Well-bred operators in $\ca_{R}$]\hfill\\
$L(k)\in\ca_{R}\otimes M_{N}(\CC)[[k]]$ is said well-bred if it 
satisfies
\bea
L_{1}(k_{1})\,a_{2}(k_2) &=& R_{21}(k_{2},k_{1})
\,a_{2}(k_2)\,L_{1}(k_{1}) \\
L_{1}(k_{1})\,a^\dag_{2}(k_{2}) &=& a^\dag_{2}(k_{2})
\,R_{12}(k_{1},k_{2})\, L_{1}(k_{1})
\eea
\end{defi}
In \cite{ZFqgp}, it has being shown that 
a well-bred operator 
$T(k)$ can be constructed as a series in $a$'s:
\be
T(k_{0}) = \II+\sum_{n=1}^\infty\, \frac{(-1)^n}{(n-1)!} 
a^\dag_{{n}\ldots {1}}\, T^{(n)}_{01\ldots n}\,a_{1\ldots n}
\label{serieT}
\ee
with 
\beano
a^\dag_{{n}\ldots {1}} &=& 
a^\dag_{\alpha_{n}}(k_{n})\ldots a^\dag_{\alpha_{1}}(k_{1})\,
e^\dag_{\alpha_{1}}\otimes\ldots e^\dag_{\alpha_{n}}\\
a_{1\ldots n} &=& a_{\beta_{1}}(k_{1})\ldots a_{\beta_{n}}(k_{n})\,
e_{\beta_{1}}\otimes\ldots e_{\beta_{n}}\\
T^{(n)}_{01\ldots n} &=& T^{(n)}_{\alpha_0,\beta_{0},\alpha_1,\beta_{1},\ldots 
\alpha_n,\beta_{n}}(k_{0},k_1,\ldots,k_n)\, 
E_{\alpha_0,\beta_{0}}\otimes 
E_{\alpha_1,\beta_{1}}\otimes\ldots\otimes E_{\alpha_n,\beta_{n}}
\nonumber
\eeano
There is an implicit summation on the indices
$\alpha_0,\beta_{0},\ldots,\alpha_n,\beta_{n}=1,\ldots,N$ and an integration over the 
spectral parameters $\int\, dk_{1}\cdots dk_{n}$. The matrices
$T^{(n)}_{01\ldots n}$ are built using only the evaluated $R$-matrix. 
For their exact expression, we refer to \cite{ZFqgp}. To clarify the 
notation, let us stress that the auxiliary spaces indices $1, 2, 
\ldots, n$ being repeated, they are "dummy" (the corresponding 
indices $\alpha$'s and $\beta$'s are summed), and as such can be 
exchanged. As a consequence, the matrices $T^{(n)}_{01\ldots n}$ 
obey the following property \cite{ZFqgp}:
\begin{prop}
Without any loss of generality,
the matrices $T^{(n)}_{01\ldots n}$ can be supposed to obey the 
following relation
\beano
\forall\ i<j && T^{(n)}_{01\ldots n}=
\cb^{-1}_{ij}\, T^{(n)}_{01\ldots n\vert ij}\,\cb_{ij}\\
\mbox{with} &&T^{(n)}_{01\ldots n}\equiv
T^{(n)}_{01\ldots n}(k_{0},k_{1},\ldots,k_{n})\\
&&T^{(n)}_{01\ldots n\vert ij}=
T^{(n)}_{0,1,\ldots,i-1,j,i+1,\ldots,j-1,i,j+1,\ldots, n}\\
\mbox{and} && \cb_{ij}=\Big(\prod_{a=i+1}^{\longleftarrow\atop 
j-1}R_{ia}(k_{i},k_{a})\Big) R_{ij}(k_{i},k_{j})
\Big(\prod_{b=i+1}^{\longrightarrow\atop j-1}
R_{bj}(k_{b},k_{j})\Big)
\eeano
Let us stress that in $T^{(n)}_{01\ldots n\vert ij}$, the spectral 
parameters $k_{i}$ and $k_{j}$ are also exchanged.
\end{prop}

The vertex operator obey the fundamental property:
\begin{prop}
The vertex operator (\ref{serieT}) obeys the quantum group 
$\cu_{R}$ relations
\be
R_{12}(k_{1},k_{2})T_{1}(k_{1})T_{2}(k_{2})=
T_{2}(k_{2})T_{1}(k_{1})R_{12}(k_{1},k_{2})
\ee
\end{prop}
Note that the construction is unique: there is only one vertex 
operator (\ie a series of the form (\ref{serieT})) which yields a well-bred 
operator.

Let us remark that the definition of well-bred operators provides 
linear equations in $L(k)$, while the quantum group relations are 
quadratic in $T(k)$. That is the reason why the use of well-bred 
operators is crucial for the construction of the quantum group.
The situation is somehow analogous to the case 
of Drinfeld twist $\cf$ for Hopf algebras: the cocycle condition is 
quadratic in $\cf$, and it is using linear equations\cite{ABRR,
millenium} in $\cf$ that 
one is able to give explicit forms for $\cf$.

\null

Due to the relation
\be
T^\dag(k)=T^{-1}(k)
\ee
the same expansion can be done for $T^{-1}(k)$:
\be
T(k_{0})^{-1} = \II+\sum_{n=1}^\infty\, \frac{(-1)^n}{(n-1)!} 
a^\dag_{{n}\ldots {1}}\, T^{(n)\dag}_{01\ldots n}\,a_{1\ldots n}
\label{Tinv}
\ee

\subsection{Integrable hierarchy associated to $\ca_{R}$}
As in the case of undeformed oscillator algebra, one can introduce 
the following Hamiltonians
\be
H_{n}=\int_{-\infty}^{\infty}dk\, k^n\, a^\dag(k)a(k),\ \ n=0,1,2,\ldots
\ee
They form an Abelian subalgebra of $\ca_{R}$, and thus define a 
hierarchy. The generators $a^\dag(k)$ and $a(k)$ are eigenvectors:
\be
{[H_{n},a^\dag(k)]}=k^n\, a^\dag(k)\mb{and}
{[H_{n},a(k)]}=-k^n\, a(k)
\ee
Indeed, as we will see in the case of the Non Linear Schr\"odinger 
equation, these generators correspond to asymptotic states of the 
physical models associated to the hierarchy.

Moreover, this integrable hierarchy is also related to the quantum 
groups $\cu_{R}$, thanks to the property:
\begin{prop}
$\cu_{R}$ generates integrals of motion for the hierarchy defined by 
the $H_{n}$'s.
\be
{[H_{n},T(k)]}=0
\ee
\end{prop}
\subsection{Fock space $\cf_{R}$}
We introduce  the Fock space $\cf_{R}$ of $\ca_{R}$, with vacuum 
$\Omega$:
\be
a(k)\Omega=0,\ \forall k\mb{;} \cf_{R}=\ca_{R}\Omega
\ee
It can be decomposed into eigenspaces of $H_{n}$:
\be
\cf_{R}=\oplus_{p=0}^\infty\int_{k_{1}\leq k_{2}\leq\ldots\leq 
k_{p}}^\oplus d^pk\, \cf_{p}(k_{1},\ldots,k_{p})
\ee
with 
\be
{[H_{n}, x]}=(k_{1}^n+k_{2}^n+\ldots+k_{p}^n)x,\ \ 
\forall x\in\cf_{p}(k_{1},\ldots,k_{p})
\ee
Moreover, since $\cu_{R}$ commutes with $H_{n}$, each 
$\cf_{p}(k_{1},\ldots,k_{p})$ eigenspace provides a representation 
of $\cu_{R}$. On this subspace, $T(k_{0})$ acts by right 
multiplication by $R_{01}(k_{0},k_{1})\ldots R_{0p}(k_{0},k_{p})$, so 
that $\cf_{p}(k_{1},\ldots,k_{p})$ can be identified with the tensor 
product of $p$ evaluation representations $\cv(k_{p})$ of $\cu_{R}$.

Note that the action of $T(k_{0})$ on $\cf_{p}(k_{1},\ldots,k_{p})$ 
allows to reconstruct the coproduct for $\cu_{R}$, although $\ca_{R}$ 
does not possess any coproduct (see \cite{ZFqgp} for more details).

\sect{Case with boundary}
\subsection{Integrable models with boundary}
We remind here the results obtained by Mintchev et al.  
\cite{Mint,bound}, 
following \cite{Cher,Skly}, on the quantum
integrable systems with boundary. As for the case without boundary, 
the central role is played by an algebra which both encodes the 
asymptotic states of the system and the effect of the boundary.
\begin{defi}[Boundary algebra $\cb_{R}$] \hfill\\
\label{defB}
    To the same $R$-matrix, one can associate another algebra, 
    the boundary algebra $\cb_{R}$, 
    with generators $\wt{a}_{i}(k)$, $\wt{a}^\dag_{i}(k)$ and $b_{ij}(k)$
    ($i,j=1,\ldots,N$) and exchange relations:
    \bea
&&\wt{a}_{1}(k_{1})\,\wt{a}_{2}(k_{2}) = R_{21}(k_{2},k_{1}) \,
\wt{a}_{2}(k_{2})\,\wt{a}_{1}(k_{1}) \label{BNl-1}\\
&&\wt{a}_{1}^\dag(k_{1})\,\wt{a}_{2}^\dag (k_{2}) = 
\wt{a}_{2}^\dag(k_{2})\,\wt{a}_{1}^\dag (k_{1}) \,R_{21}(k_{2},k_{1}) \\
&&\wt{a}_{1}(k_{1})\,\wt{a}_{2}^\dag (k_{2}) = 
\wt{a}_{2}^\dag (k_{2})\,R_{12}(k_{1},k_{2})\,\wt{a}_{1}(k_{1})+
\half{\delta_{12}}\delta(k_{1}-k_{2})+
\half{b_{12}}(k_{1})\delta(k_{1}+k_{2})\nonu
&&\wt{a}_{1}(k_{1})\,b_{2}(k_{2}) = R_{21}(k_{2},k_{1}) \, 
b_{2}(k_{2})\,R_{12}(-k_{1},k_{2})\,\wt{a}_{1}(k_{1})\\
&&b_{1}(k_{1})\,\wt{a}_{2}^\dag (k_{2}) = 
\wt{a}_{2}^\dag (k_{2})\,R_{21}(k_{2},k_{1})\,b_{1}(k_{1})\,
R_{21}(k_{2},-k_{1}) \\
&& R_{12}(k_{1},k_{2})\,\, b_{1}(k_{1})\, R_{21}(k_{2},-k_{1}) \, 
b_{2}(k_{2})\, =\,
b_{2}(k_{2})\, R_{12}k_{1},-k_{2}) \, b_{1}(k_{1})\, R_{21}(-k_{2},-k_{1}) 
\nonu
&&b(k)b(-k) = \II 
\eea
\end{defi}
We have completed the notations (\ref{not1}-\ref{not3}) by:
\bea
b_{12}(k_{1}) &=&
\sum_{i,j=1}^{N} b_{ij}(k_{1})\, e_{i}\otimes e_{j}^\dag \\
b_{1}(k_{1}) &=&
\sum_{i,j=1}^{N} b_{ij}(k_{1})\, E_{ij}\otimes\II
\mb{;}
b_{2}(k_{2})= \sum_{i,j=1}^{N} b_{ij}(k_{1})\, \II\otimes E_{ij}
\eea

The $\cb_{R}$ algebras have been introduced in \cite{bound}, where 
they were shown to 
play a fundamental role in the study of integrable systems with 
boundaries. They allow for instance the determination of off-shell 
correlation functions. As for the case without boundary, a hierarchy 
can be defined for $\cb_{R}$:
\be
\wt{H}_{2n}=\int_{0}^\infty dk\, k^{2n}\, \wt{a}^\dag(k)\wt{a}(k)
\ee
They satisfy $[H_{2n},\wt{a}^\dag(k)]=k^{2n}\,\wt{a}^\dag(k)$ and 
   $[H_{2n},\wt{a}(k)]=-k^{2n}\,\wt{a}(k)$ for $k>0$. Thus, $\wt{a}(k)$ 
   and $\wt{a}^\dag(k)$ can be regarded as asymptotic states for the 
   hierarchy with boundary.
   
\null
As for the $\ca_{R}$ algebra, one can defined an adjoint operation on 
the $\cb_{R}$ algebra:
\begin{prop}
Let $\dag$ be defined by
\be
\dag\left\{\begin{array}{lcl}
\ca_{R} & \rightarrow & \ca_{R}\\[.23ex]
{a}(k) & \mapsto & {a}^\dag(k)\\[.23ex]
{a}^\dag(k) & \mapsto & {a}(k)\\[.23ex]
{b}^\dag(k) & \mapsto & {b}(-k)
\end{array}\right.\label{Bdag}
\ee
together with $(xy)^\dag=y^\dag x^\dag$ 
$\forall\ x,y\in\cb_{R}$ and 
$R_{12}(k_{1},k_{2})^\dag=R_{21}(k_{2},k_{1})$. 

Then $\dag$ is an anti-automorphism of $\cb_{R}$.
\end{prop}

Note that there is an automorphism on $\cb_{R}$ 
 given by \cite{bound}:
\be
\rho\left\{\begin{array}{lcl}
\cb_{R} & \rightarrow & \cb_{R}\\[.23ex]
\wt{a}(k) & \mapsto & b(k)\wt{a}(-k)\\[.23ex]
\wt{a}^\dag(k) & \mapsto & \wt{a}^\dag(-k)b(-k)\\[.23ex]
b(k) & \mapsto & b(k)
\end{array}\right.\label{rho}
\ee
This automorphism encodes in algebraic terms of the 
effects of the boundary in the physical system. It is compatible with 
the adjoint operation:
\be
\rho(x^\dag)=\rho(x)^\dag,\ \ \forall\, x\in\cb_{R}
\ee

\begin{defi}[Reflection algebra $\cs_{R}$] \hfill\\
    The reflection algebra $\cs_{R}$ is the subalgebra of the boundary algebra, 
    with generators $b_{ij}(k)$
    ($i,j=1,\ldots,N$). It has exchange relations:
    \bea
&& R_{12}(k_{1},k_{2}) \, b_{1}(k_{1})\, R_{21}(k_{2},-k_{1}) \, b_{2}(k_{2})\, =\,
b_{2}(k_{2})\, R_{12}(k_{1},-k_{2}) \, b_{1}(k_{1})\, 
R_{21}(-k_{2},-k_{1})  \nonu
&&b(k)b(-k) = \II 
\eea
\end{defi}
$\cs_R$ algebras enter into the class of $ABCD$-algebras introduced 
in \cite{MF}. In the case of the nonlinear Schr\"odinger equation with 
boundary, this algebra have been introduced by Cherednik \cite{Cher}.
 They correspond, in the boundary algebra approach, 
to the symmetries of the underlying model (with boundary).
Indeed, we have the property (proved by direct calculation):
\begin{prop}
The reflection algebra generates integrals of motion of the hierarchy 
associated to the $\wt{H}_{2n}$'s:
\be
{[\cs_{R},\wt{H}_{2n}]}=0
\ee
\end{prop}
We will study in more details the reflection algebra in the following.

\subsection{Construction of $\cb_{R}$ from $\ca_{R}$\label{s:prop}}
\begin{theo}\label{bfroma}
Let $\ca_{R}$ be a ZF algebra, and $T(k)$ its corresponding well-bred 
vertex operator. Let $B(k)$ be a $N\times N$ matrix such that
    \bea
&&R_{12}(k_{1},k_{2}) \, B_{1}(k_{1})\, R_{21}(k_{2},-k_{1}) \, B_{2}(k_{2})\, =\,
B_{2}(k_{2})\, R_{12}(k_{1},-k_{2}) \, B_{1}(k_{1})\, 
R_{21}(-k_{2},-k_{1}) \nonu
&& B(k)B(-k)=\II_{N}\label{RBRB}
  \eea
Then, the following generators obey a boundary algebra $\cb_{R}^{B}$:
 \bea
\wt{a}(k) &=& \half\Big(a(k)+b(k)a(-k)\Big)\\
 \wt{a}^\dag(k) &=& \half\Big( a^\dag(k)+a^\dag(-k)b(-k)\Big)\\
 b(k) &=& T(k)B(k)T(-k)^{-1}\label{b=TT}
 \eea
 $B(k)$ is called the reflection matrix.
\end{theo}
From the form of $b(k)$ in term of $T(k)$, one can immediately deduce:
\begin{coro}
The reflection algebra $\cs_{R}$ is a subalgebra of the quantum group
$\cu_{R}$. It is also a coideal of $\cu_{R}$:
\be
\Delta b^{ab}(k)=T^{ae}(k)T^{-1}(-k)^{gb}\otimes b^{eg}(k)
\mb{\ie} \Delta \cs_{R}\subset \cu_{R}\otimes \cs_{R}
\ee
\end{coro}
\prf
From (\ref{b=TT}), one has obviously a morphism from $\cb_{R}$ 
into $\cu_{R}$. It remains to show that the kernel of this 
homomorphism is reduced to $\{0\}$. The construction has been done 
in \cite{mol} for the case of the Yangian $Y(N)$ and the corresponding 
reflection algebra. The proof in the general case just follows the 
same lines. One introduces a gradation $gr$ on $\cu_{R}$ and $\cs_{R}$
\be
gr b_{(n)}^{ab}=n\mb{and} gr T_{(n)}^{ab}=n
\ee
and shows that morphism between
the filtered algebras $gr\cs_{R}$ and $gr\cu_{R}$ has trivial kernel.
The proof relies on a PBW theorem for the quantum group $\cu_{R}$. 
For more details, see \cite{mol}. 

In the same way, still following \cite{mol}, one proves the coideal 
property. 
Using the coproduct (\ref{deltaT}), one first gets (with implicit summation 
on repeated indices):
\be
\Delta T^{-1}(k)^{ab}=T^{-1}(k)^{cb}\otimes T^{-1}(k)^{ac}
\ee
which leads to
\beano
\Delta b^{ab}(k) &=& 
\Delta\Big(T^{ac}(k)T^{-1}(-k)^{db}\Big)B^{cd}(k)\\
&=&\Big(T^{ae}(k)\otimes T^{ec}(k)\Big)
\Big(T^{-1}(-k)^{gb}\otimes T^{-1}(-k)^{dg}\Big)B^{cd}(k)\\
&=& T^{ae}(k)T^{-1}(-k)^{gb}\otimes T^{ec}(k)B^{cd}(k)T^{-1}(-k)^{dg}\\
&=& T^{ae}(k)T^{-1}(-k)^{gb}\otimes b^{eg}(k)
\eeano
\finprf

The difference between the algebras $\cb_{R}^B$ and $\cb_{R}$ can be 
seen in the following lemma:
\begin{lem}
In $\cb_{R}^B$, the automorphism $\rho$ given in (\ref{rho}) is the identity:
 \be
\wt{a}(k)=b(k)\wt{a}(-k)\mb{and}\wt{a}^\dag(k)=\wt{a}^\dag(-k)b(-k)
\label{a=ba}
\ee
\end{lem}

{\bf Remark:} Strictly speaking, $\cb_R^B$ corresponds to the coset 
(noted $\cb^\rho_R$) of the abstract boundary algebra $\cb_R$ (given by 
definition \ref{defB}) by the relation $\rho+id=0$, so that the construction of 
$\cb_R^B$ in theorem \ref{bfroma} defines 
inclusions of $\cb_R^\rho$ into $\ca_R$ (see also theorem \ref{ca/rho}).

\begin{prop}
 The Fock space $\cf_{R}$ of $\ca_{R}$ provides a Fock space representation 
 for $\cb_{R}^{B}$, defined by
\be
\wt{a}(k)\Omega=0\mb{and} b(k)\Omega=B(k)\Omega
 \ee
 \end{prop}
From the definition of the Fock space for $\ca_{R}$, one has 
$a(k)\Omega=0$, and thus $\wt{a}(k)\Omega=0$. 
The well-bred vertex operator 
$T(p)$ satisfies $T(p)\Omega=\Omega$, which implies 
$b(k)\Omega=B(k)\Omega$.

{\bf Remark:} In \cite{bound,BNLS}, the boundary algebra $\cb_{R}$ 
has several 
Fock spaces, depending on the value of $b(k)$ on $\Omega$. In the 
present article, the algebra $\cb_{R}^B$ has only one Fock space, 
but $B$ is given within the construction of $\cb_{R}^B$, and there 
are as much 
$\cb_{R}^B$-algebra constructions in the present approach, as there 
are Fock spaces in 
the approach of \cite{Cher,Skly}.

\subsection{Vertex operator construction}
One can do the same construction for the $b$ operator:
\begin{prop}[Reflection operators as vertex operators]\hfill\\
In term of the $\ca_{R}$ generators, the reflection operators read
\beano
b_{0}(k_{0}) &=& \II+\sum_{n=1}^\infty \frac{(-1)^n}{(n-1)!}
a^\dag_{{n}\ldots {1}}\, \beta^{(n)}_{01\ldots n}a_{1\ldots n}\\
\beta^{(n)}_{01\ldots n} &=& T^{(n)}_{01\ldots n}B_{0}+ 
B_{0}{T'}^{(n)\dag}_{01\ldots n}+(n-1)
\sum_{p=1}^{n-1}\left({n-2}\atop {p-1}\right) T^{(n-p)}_{0p+1\ldots n}
B_{0}{T'}^{(p)\dag}_{01\ldots p} 
\cR'_{p}
\eeano
where the prime ' indicates that one has to consider $-k_{0}$ instead 
of $k_{0}$, $B_{0}=\lim_{k\mapsto\infty}B(k)$, and
\be
\cR'_{p}=R_{p+10}'^{-1}\cdots 
R_{n0}'^{-1}=\prod_{s=p+1}^{\longrightarrow\atop 
n}R_{0s}(-k_{0},k_{s})
\ee
\end{prop}

Let us stress that the expansion is done in term of the $\ca_{R}$ 
generators $a$ and $a^\dag$, \underline{not} in term of 
$\wt{a}$ and $\wt{a}^\dag$, generators of $\cb_{R}$. It is possible 
that such an expansion would lead to a more simple expression for 
$\beta^{(n)}$.
\subsection{Hierarchy for $\cb^{B}_{R}$\label{s:hier}}
\begin{prop}[Hierarchy for $\cb_{R}^B$]\hfill\\
Let 
\be
 \wh{H}_{n}=\half\int_{-\infty}^\infty dk\, k^{n}\,\wt{a}^\dag(k)\wt{a}(k)
\ee
Then $\wh{H}_{2n+1}$ vanish identically in $\cb_{R}^B$ and
\be
\wh{H}_{2n}=\wt{H}_{2n}=\int_{0}^\infty dk\, k^{n}\,\wt{a}^\dag(k)\wt{a}(k)
\ee
   
Let $a(k)$, $a^\dag(k)$ be the generators of the $\ca_{R}$ algebra, and 
 $H_{n}$ the Hamiltonian of the 
$\ca_{R}$-hierarchy. 
Then, we have
\be
\wt{H}_{2n}={H}_{2n}+\int_{-\infty}^\infty dk\, a^\dag(k)b(k)a(k)
\ee
In other terms, one can see the Hamiltonians with boundary as the 
Hamiltonians without boundary (bulk term) plus a boundary term.
\end{prop}
\prf
The different equalities follows from the identities (\ref{a=ba}) 
 and (\ref{RBRB}).
\finprf

{\bf Remark :} The hierarchy defines integrable systems with boundary 
defined by $B(k)$. In the framework we have adopted, the definition of 
the boundary is given by the data of the 
reflection matrix $B(k)$, as it is presented in \cite{Cher}, 
\underline{but} the boundary algebra is naturally recovered here, 
contrarily to \cite{Cher}, where it is lacking for the calculation of 
off-shell correlation functions. On the other hand, in 
\cite{bound,Mint,BNLS}, the boundary 
algebra is the basic data (whence the possibility of computation of 
correlation functions), \underline{but} the data of the boundary 
condition (\ie the reflection matrix) is given with the choice of a Fock space 
$\cf_{R}^{B}$. Thus, the present framework can be viewed as a bridge 
between the approaches \cite{Cher} and \cite{bound,BNLS}.

This remark is confirmed in the following theorem (proved in 
\cite{ZFbound}):

\begin{theo}\label{ca/rho}
Let $B(k)$ be a reflection matrix of $\ca_R$, and 
$b(k)=T(k)B(k)T^{-1}(-k)$ the 
corresponding reflection operator. Let $\rho_B$ be defined by
\be
\rho_B(a(k))=b(k)\,a(-k)\mb{and}\rho_B(a^\dag(k))=a^\dag(-k)\,b(-k)
\ee
Then:

(i) $\rho_B$ is an automorphism of $\ca_R$

(ii) $\cb^B_R$ is the coset of $\ca_R$ by the ideal $Ker(\rho_B-id)$
\end{theo}

\null

We present now a property which was already proved in \cite{BNLS} for the 
case of additive $R$-matrix $R_{12}(k_{1}-k_{2})$, but the reasoning is valid 
in full generality:
\begin{prop}[Integrals of motions of the hierarchy]\hfill\\
 The reflection algebra $\cs_{R}^B$ generates integrals of motion for 
 the $\cb_{R}^B$-hierarchy.
\end{prop}

Still following the lines given in \cite{BNLS}, one gets:
\begin{prop}[Spontaneous symmetry breaking]\hfill\\
In the Fock space representation, there is a spontaneous symmetry 
breaking of the symmetry algebra through
\be
b(k)\Omega=B(k)\Omega\label{bo=Bo}
\ee
\end{prop}

\sect{Conclusion\label{conc}}
We have shown how to construct an integral hierarchy, together with 
is integrals of motion, starting from a single $R$-matrix which obeys 
to a unitary condition. It remains to construct the canonical fields 
associated to such hierarchy, as well as the physical systems which 
underlies the construction. This have been already done for the NLS 
hierarchy, and there is no doubt that the technics should apply in 
its full generality. 

As far as the case with boundary is concerned, the present technics 
has to be related to the construction of boundary states (see e.g. 
\cite{GZ,Cor}). These later relies on the existence of a reflection 
matrix. Since this matrix is also used here for the construction of 
the boundary 
algebra staring from the ZF algebra, a link between these two 
approaches should be exhibited.

An obvious generalization of this technics is the case of elliptic 
algebras, where the $R$-matrix depends also of extra parameters: if 
such a thing could be done, one would get a insight in models such as 
the XYZ model and show, as a by-product, that the elliptic 
algebras are integrals of motions of that type of models. 

Another point which need to be clarified is the correspondence with 
the studies done in \cite{delius} for Toda systems with boundaries, 
where it has been shown 
that the integrals of motions are coideals subalgebras. 


\end{document}